\documentclass[12pt,reqno]{amsart}

\newcommand{\rank}{{\operatorname{rank\,}}}

\newcommand{\sm}{\setminus}

\newcommand{\inv}{^{-1}}

\newcommand{\wt}{\widetilde}

\newcommand{\PP}{{\mathbb P}}

\newcommand{\C}{{\mathbb C}}

\newcommand{\B}{{\mathbb B}}

\newcommand{\Hom}{{\operatorname{Hom}}}
\newcommand{\codim}{{\operatorname{codim}}}

\renewcommand{\phi}{\varphi}

\newcommand{\ical}{\mathcal{I}}

\newcommand{\al}{\alpha}

\newcommand{\Ga}{\Gamma}
\newcommand{\Gr}{{\operatorname{Gr}}}

\newtheorem{theo}{{\sc Theorem}}

\newtheorem{lem}[theo]{{\sc Lemma}}

\newenvironment{rem}{\medskip\noindent{\it Remark:\/} }{\medskip}

\title[Hyperbolic hypersurfaces in $\PP^n$ of Fermat-Waring
type]{Hyperbolic hypersurfaces in $\PP^n$ of Fermat-Waring type}

\author{Bernard Shiffman}
\address{Department of Mathematics, Johns Hopkins University, Baltimore,
MD 21218, USA}
\email{shiffman@math.jhu.edu}

\author{Mikhail Zaidenberg}
\address{Universit{\'e}
Grenoble I, Institut Fourier, UMR 5582 CNRS-UJF, BP 74,
38402 St.\ Martin
d'H{\`e}res c{\'e}dex, France}
\email{zaidenbe@ujf-grenoble.fr}

\thanks{Research of the first author partially supported by NSF grant
\#DMS-9800479.}

\date{January 12, 2001}

\begin{document}

\begin{abstract} In this note we show that
there are algebraic families of hyperbolic, Fermat-Waring type hypersurfaces
in $\PP^n$ of degree
$4(n-1)^2$, for all dimensions $n\ge 2$. Moreover, there are
hyperbolic Fermat-Waring hypersurfaces  in $\PP^n$ of degree $4n^2-2n+1$
possessing complete hyperbolic, hyperbolically embedded complements.
\end{abstract}

\maketitle


Many examples have been given of hyperbolic hypersurfaces in $\PP^3$ (e.g.,
see \cite{ShZa} and the literature therein). Examples of degree 10
hyperbolic surfaces in $\PP^3$ were
recently found by Shirosaki \cite{Shr}, who also gave examples of hyperbolic
hypersurfaces with hyperbolic complements in $\PP^3$ and $\PP^4$
\cite{Shr1}. Fujimoto
\cite{Fu2} then improved Shirosaki's construction to give examples
of degree 8. Answering a question posed in
\cite{Za3}, Masuda and Noguchi \cite{MaNo} constructed the first examples of
hyperbolic projective hypersurfaces, including those with complete
hyperbolic complements,  in any dimension. Improving the degree estimates of
\cite{MaNo}, Siu and Yeung
\cite{SY} gave examples of hyperbolic hypersurfaces in
$\PP^n$ of degree
$16(n-1)^2$. (Fujimoto's recent construction \cite{Fu2} provides examples of
degree $2^n$.) We remark that it was conjectured in 1970 by S. Kobayashi
that
generic hypersurfaces in
$\PP^n$ of (presumably) degree
$2n-1$ are
hyperbolic (for $n=3$, see
\cite{DEG} and \cite{MQ}).

The following result is an
improvement of the example of Siu-Yeung
\cite{SY}: 

\begin{theo} Let  $d\ge (m-1)^2$, $m\ge 2n-1$.  Then for generic
linear functions
$h_1,\dots,h_m$ on $\C^{n+1}$, the hypersurface
$$X_{n-1} =\left\{z\in\PP^n : \sum_{j=1}^{m} h_j(z)^d=0\right\}$$
is hyperbolic.
In particular, there exist algebraic families of hyperbolic
hypersurfaces of degree $4(n-1)^2$ in $\PP^n$.
\label{maintheorem}\end{theo}

Equivalently, $X_{n-1}$ is the intersection of the Fermat hypersurface of
degree
$d$ in $\PP^{m-1}$ with a generic $n$-plane.
The low-dimensional cases $n\le 4$ of Theorem \ref{maintheorem} were given
by Shirosaki \cite{Shr1}.  Our construction is similar to those of \cite{SY}
and \cite{Shr1}.

On the other hand, examples were given in \cite{Za2} of smooth curves
of degree 5 in
$\PP^2$, with hyperbolically embedded,
complete hyperbolic complements. Examples of hyperbolically embedded
hypersurfaces in $\PP^n$ with complete hyperbolic complements were
given in \cite{MaNo} for all $n$, and in 
\cite{Shr1} for $n\le 4$ (with lower degrees). The following result
generalizes the result of Shirosaki \cite{Shr1} to all dimensions.

\begin{theo}\label{complement}  Let  $d\ge m^2-m+1$, $m\ge 2n$.  Then for
generic linear functions
$h_1,\dots,h_m$ on $\C^{n+1}$, the complement $\PP^n\sm X_{n-1}$ of the
hypersurface of Theorem~\ref{maintheorem} is complete hyperbolic and
hyperbolically embedded in $\PP^{n}$.  In particular, there exist algebraic
families of hyperbolic hypersurfaces of degree $4n^2-2n+1$ in $\PP^n$ with
hyperbolically embedded, complete hyperbolic complements.
\end{theo}

In particular, Theorem \ref{complement} provides
algebraic families of curves of degree 13 in
$\PP^2$, of surfaces of degree 31 in $\PP^3$, and so forth, whose
complements
are complete hyperbolic and hyperbolically embedded in projective space.

\medskip We shall use the following notation and lemma in our proofs of
Theorems \ref{maintheorem} and \ref{complement}:  We let
$\Gr_{m,k}$ denote the Grassmannian of complex codimension $k$ subspaces of
$\C^m$, and we write
$$Q_{m,k}=0_k\times \C^{m-k}\in \Gr_{m,k}\,.$$  Furthermore, for
$$\begin{array}{ccccccl} 1 &
\le & a &\le & c &\le & m\\ 1 & \le & b &\le & c &\le & a+b\end{array}\ \
,$$
we define

$$\Ga_{m,a,b,c}=\{V\in \Gr_{m,a}: \dim V\cap Q_{m,b} \ge m-c\}\,.$$

\begin{lem} \label{dim} $\dim \Gr_{m,a}-\dim \Ga_{m,a,b,c}= (m-c)(a+b-c)\,.$
\end{lem}

\begin{proof} Let ${\rm Mat}_{m,a}=\Hom (\C^m,\C^a)$, and let
${\rm M}_{m,a}\subset {\rm Mat}_{m,a}$
denote the surjective homomorphisms. We consider the fiber bundle
$$\begin{array}{lll} GL(a) & \longrightarrow & {\rm M}_{m,a}\\ &&\downarrow
\pi
\qquad\qquad \pi(A)=\ker A\\ && \Gr_{m,a}\end{array}\ .$$
We let
$$\wt\Ga:= \pi\inv(\Ga)= \{A\in {\rm M}_{m,a}:\dim \left(\ker
A|_{Q_{m,b}}\right)
\ge m-c\}\;;$$ whence $$\dim {\rm M}_{m,a}-\dim\wt\Ga = \dim \Gr_{m,a}-\dim
\Ga_{m,a,b,c}
\,.$$

Suppose $A\in {\rm M}_{m,a}$; i.e., $A$ is an $a\times m$ matrix of rank
$a$.  We
consider $$\wt A =\left(\begin{array}{c}I_b \quad 0 \\ \hline
A\end{array}\right) = \left(\begin{array}{cc}I_b & 0\\B &
A'\end{array}\right )\in {\rm M}_{m,a+b}\,,$$ where $B\in {\rm Mat}_{b,a},\
A'\in
{\rm Mat}_{m-b,a}$. Clearly,
$ \ker \wt A = \ker A|_{Q_{m,b}}$ and thus
$$A \in \wt \Ga\ \Leftrightarrow\ \dim (\ker \wt A)\ge m-c\ \Leftrightarrow
\ \rank \wt A \le c\ \Leftrightarrow\ \rank A' \le c-b\,.$$
It is easily seen that
$$\codim_{{\rm Mat}_{k,l}}\{C\in {\rm Mat}_{k,l}:\rank C \le r\} =
(k-r)(l-r)\,.$$
Therefore, \begin{eqnarray*} \codim_{{\rm M}_{m,a}}\wt \Ga &=&
\codim_{{\rm Mat}_{m-b,a}}\{A': \rank A' \le c-b\}
\\&=& [(m-b)-(c-b)][a-(c-b)]\\
&=&(m-c)(a+b-c)\,.\end{eqnarray*}
\end{proof}

\smallskip\noindent {\it Proof of Theorem \ref{maintheorem}:\/}
Consider the Fermat hypersurface
$$F_d:=\left\{(z_1:\dots:z_{m})\in\PP^{m-1}:\ \sum_{j=1}^{m}
z_j^d=0\right\}$$
of degree $d$ in $\PP^{m-1}$. Suppose
that $d\ge (m-1)^2,\ m\ge 2n-1$. We must show that
$X_d:=F_d\cap \PP V$ is hyperbolic for a generic
$V\in \Gr_{m,m-n-1}$.

Suppose that $f=(f_1,\dots,f_{m}):\C\to X_d$ is a holomorphic curve. By
Brody's theorem \cite{Br}, it suffices to show that $f$ is constant. We
write
$J_{m}=\{1,\dots,m\}$.  Let
$$I_0=\{j\in J_{m}:f_j=0\}\,.$$ (Of course, $I_0$ may be empty.)
We let $I_1,\dots,I_l$ denote the equivalence classes in $J_{m}\sm I_0$
under the equivalence relation
$$j\sim k \ \Leftrightarrow\ f_j/f_k = \mbox{constant}\,.$$
We let $k_{\alpha}={\rm card}\,I_{\alpha}$ and we write
$$I_\al=\{i(\al,1),\dots,i(\al,k_\al)\}\,,$$ for ${\alpha}=1,\dots,
l$, and also for $\al=0$ if $k_0\ge 1$.

The result of Toda \cite{To}, Fujimoto \cite{Fu}, and M. Green \cite{Gr}
says that  for ${\alpha}= 1,\dots, l$, we have $k_{\alpha}\ge 2$
and furthermore the constants
$$\mu_{\al j}:=f_{i(\al,j)}/f_{i(\al,1)}\in\C\sm\{0\}\qquad (1\le\al\le
l\,,\ 2\le j \le k_\al)$$ satisfy
$$1+\sum_{j=2}^{k_\al} \mu_{\al j}^d=0\,.$$
Geometrically, the image $f(\C)$ is contained in the projective $l$-plane
$Y_\mu^\ical$ given by the equations
$$z_{i(\al,j)}=\mu_{\al j}z_{i(\al,1)}\,,\quad 2\le j\le k_\al\,, 1\le
\al\le l\,; \qquad z_{i(0,j)}=0\,,\quad 1\le j\le k_0\,.$$
Here, $\ical$ denotes the partition
$\{I_0,I_1,\dots,I_l\}$ of $J_{m}$, and $\mu=\{\mu_{\al j}\}$.

Let $\wt Y^\ical_\mu\subset \C^{m}$ be the lift of $Y^\ical_\mu$.
Then
$\wt Y_\mu
^\ical\in \Gr_{m,m-l}$. If $l=1$, then $Y_\mu ^\ical$ is a point.
Otherwise, we consider $Y_\mu
^\ical\cap \PP V$ for generic $V\in
\Gr_{m,m-n-1}$.   Applying Lemma~\ref{dim} with
$$a=m-n-1,\ b=m-l,\ c=m-2$$
(changing coordinates to make $\wt
Y_\mu^\ical=Q_{m, m-l}$), we conclude that
$Y_\mu ^\ical\cap \PP V$ is either a point or is empty, i.e.  $\dim(\wt
Y_\mu^\ical\cap V)<2$, unless $V$ lies in a subvariety of $\Gr_{m,m-n-1}$ of
codimension
$$s=\big[m-(m-2)\big]\big[(m-n-1)+(m-l)-(m-2)\big]=2(m-n-l+1)\,.$$
But given a partition $\ical$, the $\mu$-moduli space of $\wt Y_\mu^\ical$
in  $\Gr_{m,m-l}$ has dimension
$$\sum_{\alpha=1}^l (k_\al-2) =m-k_0 -2l\le m-2l\,.$$

Since $m\ge 2n-1$, we have $s\ge m-2l+1$ and thus
for generic $V\in \Gr_{m,m-n-1}$, $Y_\mu
^\ical\cap \PP V$ is at most a point for all $(\ical,\mu)$. Since
$f(\C)\subset Y_\mu
^\ical\cap \PP V$ for some $Y_\mu
^\ical$, it follows that $f$ must be constant.\qed

\bigskip\noindent {\it Proof of Theorem \ref{complement}:\/} Suppose that
$d\ge m^2-m+1,\ m\ge 2n$.  Since by Theorem \ref{maintheorem}, $X_{n-1}$ is
hyperbolic for generic $h_j$, it suffices to show that any entire curve
$f: \C\to \PP^{n}\backslash X_{n-1}$ is constant (see e.g., \cite{Za2}).

We proceed as in the proof of Theorem \ref{maintheorem}.  Suppose that
$f=(f_1,\dots,f_{m}):\C\to\PP^n\sm X_{n-1}$ is a holomorphic curve. As
before, let
$I_0=\{j\in J_{m}:f_j=0\}$, and let $I_1,\dots,I_l$ denote the
equivalence classes in $J_{m}\sm I_0$ under the equivalence relation
$$j\sim k \ \Leftrightarrow\ f_j/f_k = \mbox{constant}\,.$$

Since $d>m(m-1)$, by \cite{To,Fu,Gr} we have
$$k_{\alpha}\ge 2 \ \ \mbox{and }\
1+\sum_{j=2}^{k_\al} \mu_{\al j}^d=0\ \ \mbox{for }\ 2\le\al\le l\,,$$
after permuting
the
$I_\al$ and using our previous notation. (Also, $k_1\ge 1$, but
$1+\sum_{j=2}^{k_1} \mu_{1 j}^d\ne 0$.)  The proof of this result proceeds
by
considering the map $(f_0,\dots,f_m):\C\to F_d\subset\PP^m$, where
$f_0=-\sum_{j=1}^mf^d_j=e^\phi$. (We let $I_1=\{j:f_j/f_0=
\mbox{constant}\}\ne\emptyset$.) The better estimate for
$d$ arises from the fact that $f_0$ has no zeros.

As before, the image
$f(\C)$ is contained in the projective
$l$-plane
$Y_\mu^\ical$, and
$Y_\mu ^\ical\cap \PP V$ is either a point or is empty, unless $V$ lies in a
subvariety of $\Gr_{m,m-n-1}$ of codimension $s=2(m-n-l+1)$. But this
time, the
$\mu$-moduli space of $\wt Y_\mu^\ical$ in  $\Gr_{m,m-l}$ has dimension
$$k_1-1+\sum_{\alpha=2}^l (k_\al-2) =m-k_0 -2l+1\le m-2l+1\,.$$
Since $m\ge 2n$, we have $s>  m-2l+1$ and hence for generic $V\in
\Gr_{m,m-n-1}$,
$Y_\mu ^\ical\cap \PP V$ is a point or is empty for all $(\ical,\mu)$. \qed

\begin{rem} Note that the algebraic
family of degree $d=(m-1)^2$
hyperbolic hypersurfaces in $\PP^n$ constructed in
Theorem \ref{maintheorem}
has dimension $(n+1)m-1$, as does the family
of Theorem \ref{complement}. (Recall that in Theorem \ref{maintheorem},
$m\ge
2n-1$, whereas in  Theorem \ref{complement}, $m\ge 2n$.) By the stability of
hyperbolicity theorems (see \cite{Za2}), in the corresponding
projective spaces
of degree $d$ hypersurfaces,  both families possess
open neighborhoods  consisting of hyperbolic hypersurfaces,
with hyperbolically embedded complements in the second case.
We note finally that the best possible lower bound for the
degree of a hypersurface in $\PP^n$
with hyperbolic complement should be
$d=2n+1$ (see \cite{Za1}), and the degree $2n-3$ hypersurfaces
in $\PP^n$ are definitely not hyperbolic because they contain projective
lines. (In fact, starting with $n=6$, these lines are the only rational curves
on a generic hypersurface of degree $2n-3$ in $\PP^n$; see \cite{Pa}.)
\end{rem}

\noindent {\it Acknowledgement:\/}  We would like to thank
Jean-Pierre Demailly and Junjiro Noguchi for useful discussions.  The first
author also thanks  Universit\'e Joseph Fourier, Grenoble, and  the
University of Tokyo for their hospitality.

\end{document}